\input amstex.tex
\documentstyle{amsppt}
\overfullrule=0pt
\leftheadtext{C.D.Hill and M.Nacinovich}
\rightheadtext{Elementary pseudoconcavity and fields ...}

\magnification=\magstep1

\font\sc=cmcsc10

\hyphenation{pseu-do-con-ca-ve}
\hyphenation{s-o-n-d-e-r-f-o-r-s-c-h-u-n-g-s-b-e-r-e-i-c-h}
\newcount\q
\newcount\x
\newcount\t

\def\CR{{\Cal C\Cal R}}
\long\def\se#1{\advance\q by 1
\x=0  \t=0 \bigskip
\noindent
\S\number\q \quad
{\bf {#1}}\par
\nopagebreak}

\long\def\thm#1{\advance\x by 1
\bigskip\noindent%
{\sc Theorem \number\q.\number\x}
\quad{\sl #1} \smallskip\noindent}

\long\def\thml#1#2{\advance\x by 1
\bigskip\noindent
{\sc Theorem \number\q.\number\x ({#1})}
\quad{\sl #2} \smallskip\noindent}
\long\def\lem#1{\advance\x by 1
\medskip\noindent
{\sc Lemma \number\q.\number\x}
\quad{\sl #1} \smallskip\noindent}
\long\def\prop#1{\advance\x by 1
\medskip\noindent
{\sc Proposition \number\q.\number\x}
\quad{\sl #1} \smallskip\noindent}



\long\def\form#1{\global\advance\t by 1
$${#1} \tag \number\q.\number\t$$}
\long\def\cor#1{\advance\x by 1
\bigskip\noindent%
{\sc Corollary \number\q.\number\x}%
\quad{\sl #1} \smallskip\noindent}

\def\dimo{\smallskip\noindent {\sc Proof}\quad}

\hyphenation{pseu-do-con-vex}
\hyphenation{pseu-do-con-ca-ve}
\hyphenation{Hes-si-an}
\topmatter
\title 
Elementary
Pseudoconcavity  
and 
fields of $CR$ meromorphic functions \endtitle
\author
C. Denson Hill and Mauro Nacinovich
\endauthor
\address C.D.Hill - Department of Mathematics, Stony Brook University,
Stony Brook NY 11794, USA \endaddress
\email dhill\@math.sunysb.edu \endemail
\address M.Nacinovich - Dipartimento di Matematica  -
Universit\`a di Roma ``Tor Vergata'' -
via della Ricerca Scientifica - 00133 ROMA, Italy \endaddress
\email nacinovi\@ mat.uniroma2.it\endemail
\abstract
Let $M$ be a smooth $CR$ manifold of $CR$ dimension $n$ and $CR$
codimension $k$, which is not compact, but has the local extension
property $E$. We introduce the notion of "elementary pseudoconcavity"
for $M$, which extends to $CR$ manifolds the concept of a "pseudoconcave"
complex manifold. This notion is then used to obtain generalizations,
to the noncompact case, of the results of our previous paper about
algebraic dependence, transcendence degree and related matters for
the field $\Cal K(M)$ of $CR$ meromorphic functions on $M$.
\endabstract
\endtopmatter
\document
It was A.Andreotti who, in a series of papers ([A], [AG], [ASi], 
[ASt], [AT]),
realized that a number of fundamental theorems about compact complex
manifolds (or compact analytic spaces) could be carried over to the
non-compact case.
This involved the introduction of an ``{elementary notion}'' of
pseudoconcavity: {\sl Let $X$ be a connected non-compact complex manifold.
Then $X$ is called {\it elementary pseudoconcave} if one can find
a non-empty open subset $Y\subset X$ with the following properties:
\roster
\item"($i$)" $Y$ is relatively compact in $X$, and $\partial Y$ is
smooth.
\item"($ii$)" The Levi form of $\partial Y$ restricted to the analytic
tangent space has at least one negative eigenvalue at each point of
$\partial Y$.
\endroster}
In particular, for any point $z_0\in\partial Y$ there is an analytic disc
of complex dimension $\geq 1$ which is tangent at $z_0$ to $\partial Y$,
and is contained in $Y$ except for the point $z_0$. A compact $X$ may
be considered elementary pseudoconcave, since one may take $Y=X$ and
$\partial Y=\emptyset$, so that ($i$) and ($ii$) are trivially valid.
With this condition it became possible to carry over to $X$ the basic
results of Siegel [Si] about the field $\Cal K(X)$ of meromorphic
functions on $X$. Thus one was able to discuss the transcendence degree
of $\Cal K(X)$ over $\Bbb C$, algebraic versus analytic dependence,
the fact that $\Cal K(X)$ is a simple algebraic extension of the field
$\Bbb C(f_1,\hdots,f_d)$ of rational functions of a given transcendence
basis $f_1,\hdots,f_d$, as well as an extension of Chow's theorem.\par
The present work is a sequel to our previous paper [HN12]. There we
managed to replace a compact $X$ by a {\it compact} $CR$ manifold $M$
of $CR$ dimension $n$ and $CR$ codimension $k$, which was assumed 
to have a certain local extension property $E$.
This property $E$ was taken as an axiom. In particular this axiom is
satisfied in the important case where $M$ is pseudoconcave. [Note that
here the word "{pseudoconcave}" is being used in a completely different
context, and has a completely different meaning than the 
"{elementary pseudoconcave}" mentioned above. $M$ being pseudoconcave
is a local property of $M$, having to do with 
an arbitrarily small neighborhood of each point of $M$: every complex
manifold, compact or not, is pseudoconcave in this sense.]
The purpose of the present paper is to find an analogue of the notion of
"{elementary pseudoconcavity}", for a {\sl non-compact} $CR$ manifold
$M$, and to use it to obtain generalizations of most of the results of
[HN12]. Thus we study the algebraic dependence, transcendence degree
and related matters for the field $\Cal K(M)$ of $CR$ meromorphic functions
on $M$. \par
In order to find a suitable replacement for the notion of 
"{elementary pseudoconcavity}", in the context of a $CR$ manifold
$M$ satisfying property $E$, we first discuss the complex Hessian 
$i(\partial\bar\partial)_M$ acting on real transversal $1$-jets on $M$
(see [MN]). This enables us to consider weakly and strongly pseudoconcave
domains $Y\subset M$ with smooth boundaries $\partial Y$. We prove
two lemmas which serve to replace the analytic disc, tangent to
$z_0\in\partial Y$, and otherwise contained in $Y$.\par
For more information about $CR$ manifolds, we refer the reader to
[HN1], [HN2], $\hdots$, [HN12] and to a number of interesting papers
of C.Laurent-Thi\'ebaut and J.Leiterer.\par
The first author would like to express his gratitude for the
kind hospitality of the Humboldt Universit\"at zu Berlin and
the Universit\`a di Roma "Tor Vergata".

\se{The complex Hessian on $CR$ manifolds}
In this paper $M$ will be a smooth ($\Cal C^\infty$) paracompact
manifold, whose $\roman{dim}_{\Bbb R}M=2n+k$, which has a smooth
$CR$ structure of type $(n,k)$; i.e. $n$ is the $CR-\roman{dim}_{\Bbb C}$
 and
 $k$ the $CR-\roman{codim}_{\Bbb R}$. We recall what this means: as an
abstract $CR$ manifold $M$ is really a triple $\bold M=(M,HM,J)$,
where $HM$ is a smooth real vector subbundle of rank $2n$ of the
real tangent bundle $TM$, and where $J:HM@>>>HM$ is a smooth fiber
preserving isomorphism such that $J^2=-I$. It is also required that
the {\it formal integrability conditions}
 $\left[\Gamma(M,T^{0,1}M),\Gamma(M,T^{0,1}M)\right]
\subset\Gamma(M,T^{0,1}M)$
be satisfied ($\Gamma$ means smooth sections). Here
 $T^{0,1}M=\left\{X+iJX\; \right|\; \left. X\in HM\right\}$ is the
complex subbundle of the complexification $\Bbb C HM$ of $HM$
corresponding to the eigenvalue $-i$ of $J$; we have
 $T^{1,0}M\cap T^{0,1}M= 0$ and $T^{1,0}M\oplus T^{0,1}M=\Bbb C HM$,
where $T^{1,0}M=\overline{T^{0,1}M}$. When $k=0$, we recover the
abstract definition of a {\it complex} manifold, via the
Newlander-Nirenberg theorem.\par
We denote by
 $H^0M=\{\xi\in T^*M\, | \, \langle X,\xi\rangle=0\quad\forall X\in
H_{\pi(\xi)}M\}$
the {\it characteristic bundle} of $M$. To each $\xi\in H_x^0M$,
we associate the Levi form at $\xi$:
\form{L(\xi,X)=\xi([J\tilde X,\tilde X])=d\tilde\xi(X,JX) \quad
\text{for}\quad X\in H_xM\,}
which is Hermitian for the complex structure of $H_xM$ defined by $J$.
It is of interest to carry over to  
smoothly bounded domains $Y$ of an abstract $CR$ manifold
$M$ the notions of 
pseudoconvexity and pseudoconcavity, expressed in terms of
locally defining functions,
as one has for complex
manifolds $X$, such as in the classical work of Andreotti-Grauert [AG].
Classically these notions are expressed in terms of the complex
Hessian of a locally defining function; however in the case of an abstract
$CR$ manifold the complex Hessian of a real function cannot be defined
intrinsically as a Hermitian form on $HM$. In order to obtain an invariant
notion, one must consider real transversal $1$-jets on $M$. \par
We denote by $\Cal E^j(M)$ the space of smooth complex valued exterior
forms homogeneous of degree $j$ and by 
$\Cal E^*(M)=\oplus_{j=0}^{2n+k}{\Cal E^j(M)}$
the algebra of complex valued smooth exterior forms on $M$.
Then we consider the ideal $\Cal J$ of $\Cal E^*(M)$ generated by the
one-forms vanishing on $T^{0,1}M$ and by $\overline{\Cal J}$ its
conjugate, which is the ideal of  
$\Cal E^*(M)$ generated by the
one-forms vanishing on $T^{1,0}M$.\par
The notion of transversal $1$-jet on $M$
is
best explained in terms of a choice of a splitting:
\form{\lambda:\Cal E^1(M)@>>>[\Cal J\cap\Cal E^1(M)]\oplus
[\bar{\Cal J}\cap\Cal E^1(M)]}
for the
\edef\frmba{\number\q.\number\t}
exact sequence:
\form{\CD
0\rightarrow\Cal J\cap\bar{\Cal J}\cap\Cal E^1(M) \rightarrow
[\Cal J\cap\Cal E^1(M)]\oplus
[\bar{\Cal J}\cap\Cal E^1(M)]\rightarrow\Cal E^1(M)\rightarrow0\, .
\endCD}
Note that such splittings

always exist because $\Cal J$, $\bar\Cal J$,
 $\Cal J\cap\bar\Cal J$ and $\Cal E^*(M)$  are all locally free
graded $\Cal E(M)$-modules. The splitting $\lambda=(\lambda_1,\lambda_2)$
was called a $CR$ {\it gauge} in [MN]. A $CR$ gauge is characterized by
\form{\alpha=\lambda_1(\alpha)-\lambda_2(\alpha)\, ,\qquad
\lambda_1(\alpha)\in\Cal J\cap\Cal E^1(M),\quad
\lambda_2(\alpha)\in\bar\Cal J\cap\Cal E^1(M)}
for all
\edef\frmbd{\number\q.\number\t}
$\alpha\in\Cal E^1(M)$.
\par
In a $CR$ gauge $\lambda$ a {\it real transversal $1$-jet} $\psi$
is represented by the pair $(\psi_0,\psi_1)_{\lambda}$,
where $\psi_0$ is a smooth real valued function on $M$, and $\psi_1$
is a smooth section of $H^0M$. In a different $CR$ gauge $\lambda'$,
 $\psi$ is represented by the pair $(\psi'_0,\psi'_1)_{\lambda'}$,
where
\form{\left\{\matrix\format\r&\l\\
\psi'_0&=\psi\\
\psi'_1&=\psi_1+i\lambda_1'(d\psi_0)-i\lambda_1(d\psi_0)\\
 &=\psi_1+i\lambda_2'(d\psi_0)-i\lambda_2(d\psi_0)\, .
\endmatrix\right.}
The {\it complex Hessian} $\text{{\sl i}}(\partial\bar\partial)_M\psi$ of a
real
\edef\frmbd{\number\q.\number\t}
\!{transversal} one jet $\psi$ on $M$ is then the Hermitian form on $HM$
defined by
\form{\matrix\format\r&\l\\
i(\partial\bar\partial)_M
\psi(X,JX)&=d[\psi_1-i\lambda_1(d\psi_0)](X,JX)\\
&=d[\psi_1-i\lambda_2(d\psi_0)](X,JX)\, ,
\endmatrix}
for all $X\in HM$. By (\frmbd)
\edef\frmbe{\number\q.\number\t}
the right hand side does not depend
on the choice of the $CR$ gauge $\lambda$; hence
 $i(\partial\bar\partial)_M\psi$ is an invariant notion on the abstract
 $CR$ manifold $M$. By our definition of the Levi form on $M$, we have
that
\form{L(\psi_1,X)=L(\psi_1,JX)=d\psi_1(X,JX)\, .}
Notice that if $M$ were a complex manifold ($k=0$), then we would have
\edef\frmbf{\number\q.\number\t}
$\Cal J\cap\bar\Cal J=0$, and (\frmba) would be an isomorphism, and so there
would be a natural unique choice of the $CR$ gauge; moreover $H^0M=0$,
so there is no $\psi_1$-term, and hence (\frmbe) reduces to the usual
complex Hessian of a function on the complex manifold.
Thus for an open subset $Y$ of a
$CR$ manifold $M$ with a defining function $\psi_0$,
there are many ways to extend $\psi_0$ to a real transversal
 $1$-jet $\psi$; however the various ways lead to complex Hessians
(\frmbe) which differ by a Levi form (\frmbf) of $M$. \par
Let us now consider the situation where $M$ is a generic $CR$
submanifold of a complex manifold $X$. If $f$ is a smooth real
valued function defined on $X$, we shall associate to it the real
transversal $1$-jet $\psi$ on $M$ defined by
\form{\psi=\left(f|_M,i\lambda_1(df|_M)-i(\partial f)|_M\right)|_{\lambda}
\, ,}
where $\lambda$ is some fixed $CR$ gauge on $M$. Then the complex Hessian
 $i(\partial\bar\partial)_M\psi$ of the associated transversal $1$-jet
is the restriction to $HM$ of the pullback to $M$ of the complex

Hessian $i\partial\bar\partial f$ on $X$:
\form{i(\partial\bar\partial)_M\psi=
i(\partial\bar\partial f)|_M \qquad
\text{on}\quad HM\, .}
On the other hand,
\edef\frmbh{\number\q.\number\t}
given a smooth transversal $1$-jet $\psi$ on $M$, by the Whitney
extension theorem, we can find a smooth real valued function $f$
on $X$ such that (\frmbh) holds.
\se{The $E$-property}
As in [HN12], also
in this paper we shall consider $CR$ manifolds
$M$ of type $(n,k)$ which have  property $E$ ($E$ is for
{\it extension}). We recall that
$M$ {\sl is said to have property} $E$ iff there is
an $E$-pair $(M,X)$. By an $E$-pair we mean that
\roster
\item"($i$)"\quad $M$ is a generic $CR$ submanifold of the complex
manifold $X$, and
\item"($ii$)"\quad for each $a\in M$, the restriction map induces
an isomorphism $\Cal O_{X,a} @>>>\Cal C\Cal R_{M,a}$.
\endroster
We use the notation: $\Cal O(X)$ and $\Cal O_{X,a}$
are the spaces of holomorphic functions on $X$ and
of germs of holomorphic functions defined on a neighborhood of
a point $a\in X$, respectively; likewise $\Cal C\Cal R(M)$ and 
$\Cal C\Cal R_{M,a}$ are the spaces of smooth $CR$ functions
on $M$ and the space of germs of smooth $CR$ functions 
defined on a neighborhood of
a point $a\in M$, respectively.
\smallskip
\noindent
{\sc Remark}\quad {\sl If $M$ is a pseudoconcave $CR$ manifold,
then $M$ has property $E$.} (see [HN12]).\par
\smallskip
When $k=0$, so $M$ is of type $(n,0)$, then $M$ is an $n$-dimensional  
complex manifold, and we obtain an $E$ pair by choosing $X=M$. Hence
we adopt the convention that {\it any complex manifold has property}
$E$.\par
When $n=0$, so $M$ is of type $(0,k)$, then $M$ is a smooth totally
real $k$-dimensional manifold, and we can never obtain an $E$-pair,
(unless $M=X=\text{a point}$), because then any smooth function belongs to
$\Cal C\Cal R(M)$.\smallskip
In [HN12] it was proved:
\thm{Let $(M,X)$ be an $E$-pair. Then for any open set $\omega\subset M$
there is a corresponding open set $\Omega\subset X$ such that
\roster
\item"($i$)"\quad $\Omega\cap M=\omega$, and
\item"($ii$)"\quad $r:\Cal O(\Omega)@>>>\Cal C\Cal R(\omega)$ is
an isomorphism.
\item"($iii$)" \quad 
If $f\in\Cal C\Cal R(\omega)$,
and $f$ vanishes of infinite order at $a\in\omega$, then $f\equiv 0$
in the connected component of $a$ in $\omega$.
\item"($iv$)" \quad $(r^*f)(\Omega)=f(\omega)$.
\item"($v$)"\quad If $|f|$ has a local maximum at a point $a\in\omega$, then
$f$ is constant on the connected component of $a$ in $\omega$.
\endroster}\edef\teoaa{\number\q.\number\x}
\cor{Let $(M,X)$ and $(N,Y)$ be $E$-pairs, and let $f:M@>>>N$ be a
smooth $CR$ isomorphism. Then there are $E$-pairs $(M,X')$ and
$(N,Y')$, with $X'\subset X$ and $Y'\subset Y$, such that $f$
extends to a biholomorphic diffeomorphism 
$\tilde f:X'@>>>Y'$.} 
\thm{$M$ has property $E$ if and only if for each $a\in M$, there is
an open neighborhood $\omega_a$ of $a$ in $M$ such that $\omega_a$
has property $E$.}
\se{Weakly and strongly pseudoconcave domains}
In this section we consider domains $Y\subset M$ with smooth boundaries.
We say that $Y$ is {\it weakly pseudoconcave} at
$a\in\partial Y$ if there is a smooth transversal $1$-jet
$\phi=(\phi_0,\phi_1)_{\lambda}$, defined on an open neighborhood
$U$ of $a$ in $M$ such that:
\roster
\item 
$\phi_0:U @>>>\Bbb R$ is a smooth
locally defining function for $Y$ at $a$: this means that
$d\phi_0(a)\neq 0$ and $Y\cap U=\{p\in U\, | \, \phi_0(p)<0\}$; 
\item
$i(\partial\bar\partial)_M\phi(p)+L(\xi,\,\cdot\,)$ has at least
$2$ negative eigenvalues on $H_pM$
for every $p\in U$ and every $\xi\in H^0_pM$.
\endroster

\lem{Assume that $M$ has property $E$. Let $Y$ be an open subset
of $M$ with a smooth boundary and let $U$ be a connected
relatively compact open subset of $M$.
If $Y$ is weakly pseudoconcave at every point of $\partial Y\cap 
\overline{U}$, then for every function $u$, which is defined and
$CR$ on a neighborhood of $\overline{U}$, we have:
\form{\sup_{\partial Y\cap U}{|u|}\leq \sup_{\partial (Y\cap U)
\setminus{\partial Y\cap U}}{|u|}\, .}}
\dimo
Using a partition of unity, 
we may assume that $\phi=(\phi_0,\phi_1)_{\lambda}$
is globally defined on an open neighborhood $U'$ of $\overline{U}$,
that $d\phi_0\neq 0$ on an open neighborhood $W$ of 
$\partial Y\cap\overline{U}$, and that condition (2) of the definition
is valid at each point $p\in W$.\par
We argue by contradiction. If the statement where 
false, one could find
a $CR$ function $u$, that we can assume to be defined on $U'$, and a
point $p_0\in \partial Y\cap U$ such that
\form{\mu=|u(p_0)|=\max_{\overline{Y\cap U}}{|u|}>
\sup_{\partial (Y\cap U)
\setminus{\partial Y\cap U}}{|u|}=\mu'\, .}
Set $F_\delta=\{p\in\overline{U}\, | \, \phi_0(p)\leq -\delta\}$.
Then for small $\delta>0$ the function
$\mu(\delta)=\max_{p\in F_\delta}{|u(p)|}$ is continuous
and decreasing. Choose a small $\delta\geq 0$ such that
$\mu\geq\mu(\delta)>\mu'$, 
$\partial F_\delta\setminus (\partial (Y\cap U)
\setminus{\partial Y\cap U})\subset W$, and
$\mu^2(\delta)$ is not a critical value for 
the smooth real valued function $|u|^2$, which we may do by 
Sard's theorem. Then $\mu(\delta)=|u(p_\delta)|$ for some
$p_\delta\in\partial Y_{\delta}\cap U$, because
$|u(p)|\leq \mu'<\mu(\delta)$ on $\partial{U\cap Y_\delta}\setminus
(U\cap\partial Y_\delta)\subset\partial{(U\cap Y)}\setminus
(U\cap\partial Y_\delta) $.
Then we have $d\,|u|^2(p_\delta)=k\,d\,\phi_0(p_\delta)$ with
some $k>0$, that we can arrange to be $1$.\par
Let $X$ be a complex manifold such that $(M,X)$ is an $E$ pair and
consider the open $\tilde U\subset X$ corresponding to $U$
given by Theorem 2.1. We denote by $\tilde\phi$ a smooth real valued
function on $\tilde U$ defining on $U$ the transversal $1$-jet
$\phi$. Let $\rho_1,\,\hdots,\, \rho_k$ be smooth real valued
functions on a neighborhood $\Omega$ of $p_\delta$ in $X$
such that $\partial\rho_1\wedge \cdots\wedge\partial\rho_k\neq 0$
in $\Omega$ and $M\cap\Omega=\{z\in\Omega\, | \, \rho_1=0,\,
\cdots,\,\rho_k=0 \}$. 
Then, denoting by $d_X$ the differential on $X$, for the holomorphic
extension $\tilde u$ of $u$ to $\tilde U$ we have
$d_X{|\tilde u|^2(p_\delta)}=d_X\tilde\phi(p_\delta)+
\sum_{i=1}^k{\xi_i d_X\rho_i(p_\delta)}$,
for suitable $\xi_1,\hdots,\xi_k\in\Bbb R$. Then from the inclusion:
\form{F_\delta\subset\{z\in\tilde U\, | \, |\tilde u|^2\leq\mu^2(\delta)\}}
it follows that we can find a large constant $C>0$ such that
\form{|\tilde u|^2-\mu^2(\delta)\leq \tilde\phi - \delta
+\sum_{i=1}^k{\xi_i \,\rho_i}+C\sum_{i=1}^k{\rho_i^2}}
on a neighborhood of $p_\delta$.\par
By the E.E.Levi theorem [L], all holomorphic functions defined in
$\Omega=\{z\in\tilde U\, | \,\tilde\phi - \delta
+\sum_{i=1}^k{\xi_i \,\rho_i}+C\sum_{i=1}^k{\rho_i^2}<0\}$
have a holomorphic extension to an open neighborhood of
the point $p_\delta\in\partial\Omega\cap\tilde U$. 
However, the function $z@>>>\left(\tilde u(z)-u(p_\delta)\right)^{-1}$
is holomorphic in $\Omega$ and cannot be holomorphically extended to
an open neighborhood of $p_\delta$ in $\tilde U$. This gives a contradiction,
completing the proof of the lemma.
\medskip
Next we fix a Hermitian metric $\bold h$
on the complex vector bundle
$HM@>>>M$. If $\phi=(\phi_0,\phi_1)_{\lambda}$ is a transversal $1$-jet
defined on an open subset $U$ of $M$, we can consider for each
$p\in U$ and $\xi\in H^0_pM$ the eigenvalues 
\form{\lambda_1(\phi;\xi)\leq \lambda_2(\phi;\xi)
\leq\cdots\leq\lambda_n(\phi;\xi)}
of the Hermitian form $i(\partial\bar\partial)_M\phi(p)+L(\xi,\,\cdot\,)$
with respect to $\bold h$. \par
Let $Y$ be an open subset of $M$ with smooth boundary.
We say that $Y$ is {\it strongly pseudoconcave} at
$a\in\partial Y$ if there is a smooth transversal $1$-jet
$\phi=(\phi_0,\phi_1)_{\lambda}$, defined on an open neighborhood
$U$ of $a$ in $M$ such that:
\roster
\item 
$\phi_0:U @>>>\Bbb R$ is a smooth
locally defining function for $Y$ at $a$; 
\item
there exists a positive constant $c_0>0$ 
and an open neighborhood $\omega$ of $a$ in $U$ such that 
$$\lambda_2(\phi;\xi)\leq -\, c_0<0\quad\forall p\in\omega,\;
\forall \xi\in H^0_pM\, .$$
\endroster
We obtain the following
\lem{Assume that $M$ has property $E$. Let $Y$ be an open subset
of $M$ with a smooth boundary.
If $Y$ is strongly pseudoconcave at a point $a\in \partial Y$, 
then for every open neighborhood $U$ of $a$ in $M$ there is
an open neighborhood $\omega$ of $a$ in $U$ such that
\form{\sup_{\omega}{|u|}\leq\sup_{Y\cap U}{|u|}\quad
\forall u\in\Cal C\Cal R(U)\, .}}
\dimo We can assume that the open set $U$ is so chosen that
there is a $1$-jet $\phi=(\phi_0,\phi_1)_\lambda$, 
defined on a neighborhood of $\overline{U}$, with
$\phi_0$ being a locally defining function for $Y$ at all
points of $\partial Y\cap\overline{U}$, and  
$\lambda_2(\phi;\xi)\leq -c_0<0$ for all $p\in\overline{U}$ and
$\xi\in H^0_pM$. Let $\psi_0$ be a non negative smooth real
valued function, with compact support in $U$, and with 
$\psi_0(a)>0$. Set $\phi_t=(\phi_0-t\psi,\phi_1)_\lambda$.
Then there exists $\delta>0$ such that
$\lambda_2(\phi_t;\xi)\leq -c_0/2<0$ for all $t\in\Bbb R$ with
$|t|\leq\delta$. Using the previous lemma we obtain
(\number\q.\number\t)
with $\omega=\{p\in U\, | \, \phi_0(p)-\delta\psi(p)<0\}$.
\medskip
\noindent
{\sc Example 1}\quad
Let $q\geq 3$ and consider in $\Bbb C\Bbb P^{2q-1}$ homogeneous coordinates
$(w,z)$ with $z,w\in\Bbb C^q$. Let $K$ denote the union of the two
disjoint projective $q-1$-planes:
$$K=\{w=iz\}\cup\{w=-iz\}\, .$$
Let $M$ be the $CR$ submanifold of type $(2q-3,2)$ of points
in $\Bbb C\Bbb P^{2q-1}\setminus K$ satisfying the homogeneous equations
$$\cases
\sum_{j=1}^q{|w_j|^2}=\sum_{j=1}^q{|z_j|^2}\\
\sum_{j=1}^q{w_j\,\bar z_j\, + \, z_j\,\bar w_j}=0\, .
\endcases
$$
Then, for every $1/2<\epsilon<1$ the relatively compact open subset
$$Y_\epsilon=\{\min\{|z+iw|^2,\,|z-iw|^2\} >\epsilon (|z|^2+|w|^2)\}
\cap M$$
has a smooth boundary and is strongly pseudoconcave. To see this
we observe that the boundary consists of two disjoint pieces and that
we can take the functions $\phi_\pm=\epsilon (|z|^2+|w|^2)-
|z\pm i w|^2$ as locally defining functions for $\partial Y_\epsilon$.
Note also that $M$ is $(q-2)$-pseudoconcave at each point and hence
has property $E$.

\medskip
\noindent
{\sc Example 2}\quad
Denote by $w_0,w_1,w_2,z_0,z_1,z_2$ the homogeneous coordinates
of $\Bbb C\Bbb P^5$. Let $Q$ be the 
ruled quadric described by the homogeneous equations
$$w_0=0,\; z_0=0,\; w_1z_2=w_2z_1\, .$$
Then consider the noncompact $CR$ manifold $M$ consisting of
the points of $\Bbb C\Bbb P^5\setminus Q$ whose homogeneous coordinates
satisfy:
$$\cases
w_0\bar w_1+w_1\bar w_0=z_0\bar z_1+z_1\bar z_0\\
w_0\bar w_2+w_2\bar w_0=z_0\bar z_2+z_2\bar z_0
\endcases
$$
Then $M$ is $1$-pseudoconcave
at all points which do not
belong to the $3$-plane
$\Sigma=\{w_0=0,\; z_0=0\}$. The Levi form is identically zero at the
points of $\Sigma\cap M$, because $\Sigma\setminus Q\subset M$.
Let $d$ be the distance in the Fubini-Study
metric of $\Bbb C\Bbb P^5$ and denote by $Y_\epsilon$ the set of points
of $M$ having a distance $>\epsilon$ from $Q$. Since $Q$ is smooth,
for small $\epsilon>0$ the boundary of $Y_\epsilon$ is smooth and one
can verify, by using the function $\phi(p)=\epsilon-d(p,Q)$, that
for small $\epsilon>0$ the domain $Y_\epsilon$ is strongly pseudoconcave
at all points $a\in\partial Y_\epsilon$.
Note that $M$ does not have the property $E$, because it is not minimal
at the points of $\Sigma\setminus Q$.

\se{$CR$ meromorphic functions on 
elementary pseudoconcave $CR$ manifolds.}
Let $M$ be a connected smooth $CR$ manifold of type $(n,k)$.
We say that $M$ is {\it elementary pseudoconcave} if it contains
a relatively compact non empty open subset $Y$,
with a smooth boundary which is
strongly pseudoconcave at every point $a\in\partial Y$.
\par
Note that for a compact $M$ we can take $Y=M$, so that $\partial Y=\emptyset$
and the condition above is trivially satisfied: hence compact $CR$
manifolds are trivially elementary pseudoconcave. The $CR$
manifolds $M$ described in the Examples 1,2 
at the end of the previous section
provide  examples of noncompact elementary pseudoconcave
$CR$ manifolds, the first having and the second not having
property $E$.\par
For the notion of $CR$ meromorphic functions we refer to our previous
article [HN12].
\thm{Let $M$ be a connected smooth $CR$ manifold of type $(n,k)$.
If $M$ has property $E$ and is elementary pseudoconcave, then
the field $\Cal K(M)$ of $CR$ meromorphic functions on $M$ has
transcendence degree over $\Bbb C$ 
less or equal to $n+k$.}
\par
\quad Setting $k=0$ above, we recover Andreotti's
generalization [A] of Satz 1 of Siegel [Si].
\dimo
The statement means: Given
$n+k+1$ $CR$ meromorphic functions \par
\centerline{$f_0,\; f_1, \;\hdots,\; f_{n+k}$\quad on $M$,}\par\noindent
there exists a non zero polynomial with complex coefficients
$F(x_0,x_1,\hdots,x_{n+k})$ such that
\form{F(f_0,f_1,\hdots,f_{n+k})\equiv 0\quad\text{on}
\quad M\, .}\edef\formba{\number\q.\number\t}
Because of property (E) we may regard $M$ as a generic $CR$
submanifold of an $n+k$ dimensional complex manifold $X$.\par
For each point $a\in M$ there is a connected open coordinate 
neighborhood $\Omega_a$, in which the holomorphic coordinate
$z_a$ is centered at $a$. We choose $\Omega_a$ in such a way that
$\omega_a=\Omega_a\cap M$ is a connected neighborhood of $a$ in $M$.
Moreover we can arrange that, for $j=0,1,\hdots,n+k$, each $f_j$
has a representation
\form{f_j=\dsize\frac{p_{ja}}{q_{ja}}\quad\text{on}\quad\omega_a}
with 
$p_{ja}$ and $q_{ja}$ being smooth $CR$ functions in $\omega_a$.
According to Theorem \teoaa\; we may also assume that the restriction 
map $\Cal O(\Omega_a)@>>>\CR(\omega_a)$ is an isomorphism.
For each $CR$ function $g$ on $\omega_a$, we denote its unique
holomorphic extension 
to $\Omega_a$ by $\tilde{g}$.
By a careful choice of the $p_{ja}$ and $q_{ja}$, and an additional
shrinking of $\omega_a$, $\Omega_a$, we can also arrange that
\form{\tilde{f}_j=
\dsize\frac{\tilde{p}_{ja}}{\tilde{q}_{ja}}\quad\text{on}
\quad\Omega_a\, ,}
with the functions $\tilde{p}_{ja}$ and $\tilde{q}_{ja}$
being holomorphic and having no 
nontrivial common factor at each point in a neighborhood of
$\overline{\Omega}_a$. For each pair of points $a,b$ on $M$ we have the
transition functions
\form{\tilde{q}_{ja}=g_{jab}\tilde{q}_{jb}\, 
,}\noindent
which are holomorphic and non vanishing on a neighborhood of
$\overline{\Omega}_a\cap\overline{\Omega}_b$. 
\par
Let $Y$ be a relatively compact open subset of $M$ with $\partial Y$
smooth 
and strongly pseudoconcave at each point.
For each $a\in \overline{Y}$ 
we can choose the polydiscs:
\form{K_a=\{|z_a|\leq r_a\}\quad\text{and}\quad L_a=\{|z_a|<e^{-1}r_a\}\, 
,}\noindent
where $|z_a|$ denotes the max norm in $\Bbb C^{n+k}$, and $r_a>0$,
so that $K_a\Subset\Omega_a$
and for all $u\in\CR(\omega_a)$ we have
\form{\sup_{K_a}{|\tilde u|}\leq\sup_{Y\cap\omega_a}{|u|}\, .}
This is trivial when $a\in Y$, as we can take 
$K_a\subset\widetilde{Y\cap\omega_a}$
in this case. When $a\in\partial Y$, we apply
Lemma 3.2 to find an open neighborhood $\omega$ of $a$ in $\omega_a$
so  that (3.6) is valid with $U=\omega_a$, and next we choose
$K_a\subset\Omega$, where
$\Omega$ is the open set in $X$ of Theorem 2.1.\par
By the compactness of $\overline{Y}$, we may
fix a finite number of points $a_1,a_2,\hdots,a_m$ on $M$, such that the
$L_{a_1},L_{a_2},\hdots,L_{a_m}$ provide an open 
covering of $\overline{Y}$. 
Then
we choose positive real numbers $\mu$ and $\nu$ to provide the bounds:
\form{|g_{0ab}|<e^\mu\qquad\text{and}\qquad \left|\dsize\prod_{j=1}^{n+k}{
g_{jab}}\right|<e^\nu} \edef\formbf{\number\q.\number\t}
on $\overline{\Omega}_a\cap\overline{\Omega}_b$ for
$a,b=a_1,a_2,\hdots,a_m$.\par
Consider a polynomial with complex coefficients to be determined later,
\hfill\linebreak
$F(x_0,x_1,\hdots,x_{n+k})$ of degree $s$ with respect to $x_0$ and
of degree $t$ with respect to each $x_i$ for $i=1,2,\hdots,n+k$. The number
of coefficients to be determined is
\form{A=(s+1)\cdot (t+1)^{n+k}\, .}
Now, letting $a$ stand for any one of the $a_1,a_2,\hdots,a_m$, we
introduce the functions
\form{Q_a=\tilde q^s_{0a}\prod_{j=1}^{n+k}{\tilde q_{ja}^t}\,,
\quad P_a=Q_aF(\tilde f_0,\tilde f_1,
\hdots,\tilde f_{n+k})}\noindent
which are holomorphic on a neighborhood of $\overline{\Omega}_a$.
For a positive integer $h$, to be made precise later, we wish to impose
the condition, for $a=a_1,a_2,\hdots,a_m$, that $P_a$ vanishes to order
$h$ at $a$. In terms of our local coordinates $z_a$, this means that
all partial derivatives of order $\leq h-1$ must vanish at $z_a=0$.
This imposes a certain number of linear homogeneous conditions on the
unknown coefficients of the polynomial $F$. The number of such conditions
is 
\form{B=m\binom{n+k+h-1}{n+k}\leq m\, h^{n+k}\, .}
If we can arrange that $B<A$, then this system of linear homogeneous
equations has a non trivial solution.
\par
However, in order to apply the Schwarz lemma later, we need also to
arrange that $s$, $t$ and $h$ satisfy
\form{\mu\,s\,+\,\nu\,t\, <\, h\, .}\edef\formbl{\number\q.\number\t}
To this end we fix $s$ to be an integer with $s>m\nu^{n+k}$. Thus, for 
each positive $h$, we denote by $t_h$ the largest positive integer 
satisfying
$st_h^{n+k}<mh^{n+k}$.
In this way we obtain that
\form{B\leq mh^{n+k}\leq s\left(t_h+1\right)^{n+k}<(s+1)
\left(t_h+1\right)^{n+k}=A\, .}
On the other hand, since $t_h@>>>\infty$ as $h@>>>\infty$, by choosing
$h$ sufficiently large we have
\form{m\left(\dsize\frac{\mu s}{t_h}+\nu\right)^{n+k}
<s\, ,}
which implies (\formbl) for $t=t_h$. Set
\form{\Upsilon=\max_{1\leq i\leq m}\,
\max_{K_{a_i}}{|P_{a_i}|}\, .}
This maximum is obtained at some point $z^*$ belonging to some
$K_{a^*}$, for $a^*$ equal to some one of $a_1,a_2,\hdots,a_m$.
Since $z^*\in K_{a^*}\subset\Omega_{a^*}$, because of our choices of
the $\omega_a$, $\Omega_a$, according to ($iv$) in Theorem 2.1
and (4.6),
there is another point $z^{**}\in \omega_{a^*}\cap\overline{Y}$ 
such that
\form{P_{a^*}(z^*)=P_{a^*}(z^{**})\, .}\edef\formbp{\number\q.\number\t}
But the point $z^{**}$ belongs to some $L_{a^{**}}\subset K_{a^{**}}$,
where $a^{**}$ is one of the $a_1$, $a_2$, $\hdots$,
$a_m$. Hence by the Schwartz
lemma of Siegel [Si] we obtain
\form{\left|P_{a^{**}}(z^{**})\right|
\leq \Upsilon\, e^{-h}\, .} \edef\formbq{\number\q.\number\t}
However
\form{P_{a^{*}}(z^{**})=P_{a^{**}}(z^{**})\left[
g_{0a^*a^{**}}^s(z^{**})\dsize\prod_{j=1}^{n+k}{g^t_{ja^*a^{**}}(z^{**})}
\right]\, .} 
Hence from (\formbf), (\formbp), (\formbq) we obtain
\form{\Upsilon=\left|P_{a^{*}}(z^{**})\right|\leq \Upsilon\, 
e^{\mu s+\nu t -h}\, .}
By (\formbl) this implies that $\Upsilon=0$. Hence each $P_{a_j}\equiv 0$,
which in turn yields $F(\tilde f_0,\tilde f_1,\hdots,\tilde f_{n+k})\equiv 0$.
Therefore restricting to $M$ we get (\formba). This completes the proof.
\medskip
We note that our proof follows closely that of the corresponding result
(Theorem 2.1) of [HN12], the only change consisting in the restriction of
the covering by  polycylinders
to the points of the closure of the subdomain $Y$ of $M$, and the use
of Lemma 3.2 to reduce again the discussion to polycylinders centered
at points of $\overline{Y}$.
\par
In a completely similar way we can extend also the other results of
[HN12] to the present situation. We shall therefore refer, for the
proofs of the following results, to [HN12], as only small changes,
similar to those explained in 
detail in the proof of Theorem \number\q.\number\x
\; are needed.
\smallskip
Let $f_0,f_1,\hdots,f_\ell\in \Cal K(M)$. We recall that they are
{\it analytically dependent} if
\form{df_0\wedge df_1\wedge \cdots \wedge df_\ell=0\quad\text{
where it is defined.}}
\thm{Let $M$ be a connected smooth 
elementary pseudoconcave $CR$ manifold of type
$(n,k)$, having property $E$. Let $f_0,f_1,\hdots,f_\ell\in \Cal K(M)$.
Then they are algebraically dependent over $\Bbb C$ if and only if they
are analytically dependent.}

\thm{Let $M$ be a connected smooth 
elementary pseudoconcave $CR$ manifold of type
$(n,k)$, having property $E$.
Let $d$ be the transcendence degree of  $\Cal K(M)$ over $\Bbb C$,
and let $f_1,f_2,\hdots, f_d$ be a maximal set of algebraically independent
$CR$ meromorphic functions in $\Cal K(M)$. Then $\Cal K(M)$ is a
simple finite algebraic extension of the field
$\Bbb C(f_1,f_2,\hdots, f_d)$ of rational functions of 
$f_1,f_2,\hdots, f_d$.}
Setting $k=0$ above, and taking the special case where $d=n$, we recover
Andreotti's generalization [A] of
Satz 2 of Siegel [Si].\par
As in [HN12], this theorem can be 
derived from the following:
\prop{Let $f_1,f_2,\hdots,f_\ell$ be $CR$ meromorphic functions in 
$\Cal K(M)$. Then there exists a positive integer 
$\kappa=\kappa(f_1,f_2,\hdots,f_\ell)$ such that every $f_0\in\Cal K(M)$,
which is algebraically dependent on $f_1,f_2,\hdots,f_\ell$, satisfies a
nontrivial polynomial equation of degree $\leq\kappa$ whose coefficients
are rational functions of $f_1,f_2,\hdots,f_\ell$.}
\edef\teordb{\number\q.\number\x}
\smallskip
In particular, fix a
maximal set $f_1,f_2,\hdots, f_d$ of algebraically independent
$CR$ meromorphic functions on $M$, where $d$ is the transcendence
degree
of $\Cal K(M)$. Consider an $f\in\Cal K(M)$. Then $f$ is algebraically
dependent on $f_1,f_2,\hdots, f_d$; i.e. it satisfies an equation
\form{f^\lambda+g_1f^{\lambda-1}+\cdots+g_\lambda=0\, ,
}
\noindent
where $g_1,g_2,\hdots,g_\lambda\in\Bbb C(f_1,f_2,\hdots, f_d)$.
The minimal $\lambda$ for which such an equation holds is called the
{\it degree} of $f$ over $\Bbb C(f_1,f_2,\hdots, f_d)$. 
By Proposition \teordb\; this degree is bounded from above by
$\kappa=\kappa(f_1,f_2,\hdots, f_d)$. Now choose an element
$\Theta\in\Cal K(M)$ so that its degree $\alpha$ is maximal. For any
$f\in\Cal K(M)$ consider the algebraic extension field
$\Bbb C(f_1,f_2,\hdots, f_d,\Theta,f)$. By the primitive element
theorem this extension is simple; i.e. there exists an element
$h\in\Bbb C(f_1,f_2,\hdots, f_d,\Theta,f)$ such that
$\Bbb C(f_1,f_2,\hdots, f_d,\Theta,f)= \Bbb C(f_1,f_2,\hdots, f_d,h)$.
Then 
\form{\matrix\format\r&\l\\
 \alpha&\qquad \geq \qquad
[\Bbb C(f_1,f_2,\hdots, f_d,h):\Bbb C(f_1,f_2,\hdots, f_d)]\\
\\
&=[\Bbb C(f_1,f_2,\hdots, f_d,\Theta,f):\Bbb C(f_1,f_2,\hdots, f_d,\Theta)]\\
&\qquad\qquad\qquad \times \,
[\Bbb C(f_1,f_2,\hdots, f_d,\Theta):\Bbb C(f_1,f_2,\hdots, f_d)]\\
\\
&\qquad \qquad \qquad\qquad \geq\alpha\, .\\
\endmatrix}
Hence the first factor on the right must be one; therefore
$f\in\Bbb C(f_1,f_2,\hdots, f_d,\Theta)$. The conclusion is that
\form{\Cal K(M)= \Bbb C(f_1,f_2,\hdots, f_d,\Theta)=
\Bbb C(f_1,f_2,\hdots, f_d)[\Theta]\, ,
}
\noindent
and {\sl any $f\in\Cal K(M)$ can be written as a polynomial of degree
$<\alpha$ having coefficients that are rational functions of
$f_1,f_2,\hdots, f_d$.}\par
From the above remark we derive the
\prop{There is an open neighborhood $U$ of $M$ in $X$ such that the
restriction map
\form{ \Cal K(U) @>>>\Cal K(M)}
\noindent
is an isomorphism.}
Here $\Cal K(U)$ denotes the field of meromorphic functions on $U$.
\smallskip
Let $M$ be a connected smooth abstract $CR$ manifold of type $(n,k)$.
Consider a complex $CR$ line bundle $F@>\pi>>M$ over $M$. Introduce
the graded ring
\form{\Cal A(M,F)=\bigcup_{\ell=0}^\infty{\CR(M,F^\ell)}\, ,
}
\noindent
where $\CR(M,F^\ell)$ are the smooth global $CR$ sections of the
$\ell$-th tensor power of $F$. Note that if $\sigma_1\in\CR(M,F^{\ell_1})$
and $\sigma_2\in\CR(M,F^{\ell_2})$, then
$\sigma_1\sigma_2\in\CR(M,F^{\ell_1+\ell_2})$.
\par
Assume that we are in a situation where smooth sections of $F$ have the weak
unique continuation property; e.g. we could take $M$ to be
essentially pseudoconcave (see [HN8]). Then
$\Cal A(M,F)$ is an integral domain because $M$ is connected. Let
\form{\Cal Q(M,F)=\left\{\dsize\frac{\sigma_1}{\sigma_2}\, \left|
\, \sigma_1,\sigma_2\in\CR(M,F^\ell)\;\text{for some $\ell$, and}\;
\sigma_0\not\equiv 0\right.\right\}}
\noindent
denote the field of quotients. \par
Then
\footnote{See [HN12]:We can associate a $CR$ meromorphic
function $f$ to any pair $(p,q)$, where $p$ and $q$ are smooth global
$CR$ sections of a smooth complex $CR$ line bundle $F@>\pi>>M$, with
$q\not\equiv 0$. Another pair $(p',q')$, which are smooth
$CR$ global sections of another such $F'@>{\pi'}>>M$, with $q'\not\equiv 0$,
define the same $f$ iff $pq'=p'q$ as sections of $F\otimes F'$.
Note that $f=p/q$ is a well defined smooth $CR$ function where $q\neq 0$.
With this more general definition, we get a new collection
$\hat{\Cal K}(M)$ of objects called $CR$ meromorphic functions on $M$.
Observe that $\hat{\Cal K}(M)$ is a field. For an essentially
pseudoconcave $M$, which has property $E$, $\Cal K(M)$ is a subfield
of $\hat{\Cal K}(M)$. If in addition $M$ is $2$-pseudoconcave, then all
smooth complex $CR$ line bundles over $M$ are locally $CR$ trivializable,
and then $\Cal K(M)=\hat{\Cal K}(M)$.}
$\Cal Q(M,F)\subset\hat{\Cal K}(M)$,
and $\CR(M)=\Cal A(M,\text{trivial bundle})$.
\prop{Assume that $M$ is elementary pseudoconcave and has property $E$. 
\roster
\item 
If $F$ is locally
$CR$ trivializable, then $\Cal Q(M,F)$ is an algebraically closed subfield
of $\Cal K(M)$.
\item 
There exists a choice of a locally $CR$ trivializable $F$ such that
$\Cal Q(M,F)=\Cal K(M)$.
\par
Assume moreover that $M$ is essentially pseudoconcave
\footnote{See [HN8]: $M$ is {\it essentially pseudoconcave} 
iff it is {\it minimal}, i.e. does not
contain germs of $CR$ manifolds with the same $CR$ dimension and
a smaller $CR$ codimension,
and admits a Hermitian
metric on $HM$ for which the traces of the Levi forms are zero at each
point.}. Then
\item $\Cal Q(M,F)$ is algebraically closed in $\hat{\Cal K}(M)$.
\par
In case $M$ is compact and satisfies both hypothesis, then
\item ${\Cal K}(M)$ is algebraically closed in $\hat{\Cal K}(M)$.
\endroster}
\smallskip
Let $M$ be a connected smooth 
elementary pseudoconcave $CR$ manifold of type $(n,k)$,
having property $E$. Then
\thm{Let $F@>\pi>> M$ be a locally $\CR$ trivializable smooth
complex $CR$ line bundle over $M$. Then
\form{\roman{dim}_{\Bbb C}{\CR(M,F)}<\infty\, 
.}\noindent}
\thm{Let $M$ be a paracompact connected elementary pseudoconcave
$CR$ manifold of type $(n,k)$ having property $E$.
Let $\tau:M@>>>\Bbb C\Bbb P^N$ be a smooth $CR$ map. Suppose that
$\tau$ has maximal rank $2n+k$ at one point of $M$. Then
$\tau(M)$ is contained in an irreducible algebraic subvariety
of complex dimension $n+k$, and the transcendence degree of
$\Cal K(M)$ over $\Bbb C$ is $n+k$.}

\thm{The following are equivalent:
\roster
\item
$M$ has a smooth $CR$ embedding as a locally closed
$CR$ submanifold of some 
$\Bbb C\Bbb P^N$.
\item There exists over $M$ a smooth complex $CR$ line bundle $F$ such that
the graded ring $\Cal A(M,F)=\bigcup_{\ell=0}^\infty{\CR(M,F^\ell)}$
separates points and gives ``local coordinates'' at each point of $M$.
\endroster}


\enddocument


\Refs\widestnumber\key{aaaaa}

\hyphenation{con-ca-ve}
\ref\key A
\by A.Andreotti
\paper Th\'eor\`emes de dependence alg\'ebrique sur les espaces
complexes pseudo-concaves
\jour Bull. Soc. Math. France
\vol 91
\yr 1963
\pages 1-38
\endref

\ref\key ASi
\by A.Andreotti, Y-T.Siu
\paper Projective embedding of pseudoconcave spaces
\jour Ann. Scuola Norm. Sup. Pisa
\yr 1970
\vol 24 (s. 3)
\pages 231-178
\endref

\ref\key ASt
\by A.Andreotti, W.Stoll
\paper Analytic and algebraic dependence of meromorphic functions
\jour Lecture Notes in Mathematics
\yr 1971
\vol 234
\publ Springer
\publaddr Berlin
\endref

\ref\key AT
\by A.Andreotti, G.Tomassini
\paper Some remarks on pseudoconcave manifolds 
\inbook Essays on topology and related topics
\bookinfo M\'emoirs dedi\'es \`a G. deRham
\pages 85-104
\yr 1970
\endref



\ref\key AG
\by A.Andreotti, H.Grauert
\paper Algebraische K\"orper von automorphen Funktionen
\jour Nachr. Ak. Wiss. G\"ottingen
\yr 1961
\pages 39-48
\endref



\ref\key BHN
\by J.Brinkschulte, C.D.Hill, M.Nacinovich
\paper Remarks on weakly pseudoconvex boundaries
\pages 1-10
\yr 2003
\vol 14
\jour Indag. Mathem. N.S.
\endref

\ref\key BP
\by A.Boggess, J.Polking
\paper Holomorphic extensions of $CR$ functions
\jour Duke Math. J.
\vol 49
\yr 1982
\pages 757-784
\endref

\ref\key C
\by W.L.Chow
\paper On complex compact analytic varieties
\jour Amer. J. Math
\vol 71
\yr 1949
\pages 893-914
\endref


\ref\key DCN
\by L.De Carli,
M.Nacinovich
\paper Unique continuation in abstract $CR$ manifolds
\jour
Ann. Scuola Norm. Sup. Pisa Cl. Sci.
\vol 27
\yr 1999
\pages
27-46
\endref

\ref \key HN1
\by C.D.Hill, M.Nacinovich
\paper A necessary condition for global Stein immersion of compact
$CR$ manifolds
\jour Riv. Mat. Univ. Parma
\vol 5
\yr
1992
\pages 175-182
\endref

\ref \key HN2
\bysame
\paper The topology of Stein CR manifolds and the Lefschetz theorem
\jour Ann. Inst.
Fourier, Grenoble
\vol 43
\yr 1993
\pages 459--468
\endref


\ref \key HN3
\bysame 
\paper Pseudoconcave $CR$ manifolds
\inbook Complex Analysis and
Geometry, (eds Ancona, Ballico, Silva)
\publ Marcel Dekker, Inc
\yr
1996
\publaddr New York
\pages 275--297
\endref


\ref \key HN4
\bysame
\paper Aneurysms of pseudoconcave $CR$ manifolds
\jour Math.
Z.
\vol 220
\pages 347--367
\yr 1995
\endref

\ref \key HN5
\bysame
\paper
Duality and distribution cohomology of $CR$ manifolds
\jour Ann. Scuola
Norm. Sup. Pisa
\vol 22
\yr 1995
\pages 315--339
\endref


\ref \key HN6
\bysame
\paper On the Cauchy problem in complex analysis
\jour Annali
di matematica pura e applicata
\vol CLXXI (IV)
\yr 1996
\pages
159-179
\endref


\ref\key HN7
\bysame
\paper Conormal
suspensions of differential complexes
\jour J. Geom. Anal. 
\vol 10
\yr
2000
\pages 481-523
\endref

\hyphenation{qua-der-ni}
\ref\key HN8
\bysame
\paper A weak
pseudoconcavity condition for abstract almost $CR$ manifolds
\jour Invent.
math.
\yr 2000
\vol 142
\pages 251-283
\endref

\ref\key HN9
\bysame
\paper Weak pseudoconcavity and the maximum modulus principle
\jour Annali di Matematica 
\vol 182
\pages 103-112
\yr 2003
\endref

\ref\key HN10
\bysame
\paper Pseudoconvexity at infinity
\inbook Selected topics in Cauchy-Riemann geometry (ed. S.Dragomir)
\yr  2001
\vol 9
\pages 139-173
\publ Aracne
\publaddr Caserta
\endref


\ref\key HN11
\bysame
\paper Two lemmas on double complexes and their applications to
$CR$ cohomology
\inbook Selected topics in Cauchy-Riemann geometry (ed. S.Dragomir)
\yr  2001
\vol 9
\pages 125-138
\publ Aracne
\publaddr Caserta
\endref

\ref\key HN12
\bysame
\paper Fields of $CR$ meromorphic functions
\jour Sem. Mat. Univ. Padova
\vol 111
\yr 2004
\pages
\endref

\ref\key HP
\by C.D.Hill, E.Porten
\paper The H-principle and pseudoconcave $CR$ manifolds
\inbook Proceedings of Albertofest
\publaddr Cuernavaca, Mexico
\yr 2004
\endref


\ref \key L
\by E.E.Levi
\paper Studii sui punti singolari essenziali delle funzioni analitiche
di due o pi\`u variabili complesse
\jour Ann. Mat. Pura Appl.
\vol XVII (s.III)
\yr 1909
\moreref Opere, Cremonese, Roma, 1958, 187-213
\endref

\ref \key MN
\by C.Medori, M.Nacinovich
\paper Pluriharmonic functions on abstract $CR$ manifolds
\jour Annali di Matematica pura ed applicata
\vol CLXX (s.IV)
\yr 1996
\pages 377-394
\endref




\ref \key NV
\by
M.Nacinovich, G.Valli
\paper Tangential Cauchy-Riemann complexes on
distributions
\jour Ann. Mat. Pura Appl.
\vol 146
\yr 1987
\pages
123-160
\endref

\ref \key Se
\by J.P.Serre
\paper Fonctions automorphes, quelques majorations dans le cas o\`u
$X/G$ est compact.
\inbook S\'eminaire H. Cartan 1953-54
\yr 1957
\publ Benjamin
\publaddr New York
\endref

\ref\key Si
\by C.L.Siegel
\paper Meromorphe Funktionen auf kompakten analytischen
Manningfaltigkeiten
\jour Nachr. Ak. Wiss G\"ottingen
\yr 1955
\pages 71-77
\endref
\endRefs
\enddocument